\documentclass[12pt,draftcls,onecolumn]{IEEEtran}
\input epsf

\usepackage{amsmath,amssymb}

\usepackage{color}

\newcommand{\col}{\mbox{col}}

\newcommand{\nnum}{\nonumber}
\newcommand{\EQ}{\begin{eqnarray}}
\newcommand{\EN}{\end{eqnarray}}
\newcommand{\EQQ}{\begin{eqnarray*}}
\newcommand{\ENN}{\end{eqnarray*}}

\newcommand{\bremark}{\begin{remark} \begin{rm} }
\newcommand{\eremark}{ \end{rm} \rule{1mm}{2mm}
\end{remark} }
\newcommand{\btheorem}{\begin{theorem} \begin{rm} }
\newcommand{\etheorem}{ \end{rm} \rule{1mm}{2mm}
\end{theorem} }
\newcommand{\blemma}{\begin{lemma} \begin{rm} }
\newcommand{\elemma}{ \end{rm} \rule{1mm}{2mm}
\end{lemma} }
\newcommand{\bcorollary}{\begin{corollary} \begin{rm} }
\newcommand{\ecorollary}{ \end{rm} \rule{1mm}{2mm}
\end{corollary} }
\newcommand{\bdefinition}{\begin{definition}\begin{rm} }
\newcommand{\edefinition}{ \end{rm} \rule{1mm}{2mm}
\end{definition} }
\newcommand{\bproposition}{\begin{proposition} \begin{rm} }
\newcommand{\eproposition}{ \end{rm} \rule{1mm}{2mm}
\end{proposition} }
\newcommand{\bexample}{\begin{example} \begin{rm} }
\newcommand{\eexample}{ \end{rm} \rule{1mm}{2mm}
\end{example} }
\newcommand{\basm}{\begin{assumption} \begin{rm}}
\newcommand{\easm}{\end{rm} 
\end{assumption}}

\newcommand{\real}{\mathbb{R}}

\newcommand{\SF}{\operatorname{SF}}

\newcommand{\AF}{\operatorname{AF}}
\newcommand{\AS}{\operatorname{AS}}
\newcommand{\AR}{\operatorname{AR}}

\newtheorem{theorem}{\bf Theorem}[section]
\newtheorem{lemma}{\bf Lemma}[section]
\newtheorem{definition}{\bf Definition}[section]
\newtheorem{remark}{\bf Remark}[section]
\newtheorem{corollary}{\bf Corollary}[section]
\newtheorem{proposition}{\bf Proposition}[section]
\newtheorem{example}{\bf Example}[section]
\newtheorem{assumption}{\bf Assumption}[section]

\newcommand\oprocendsymbol{\hbox{$\bullet$}}
\newcommand\oprocend{\relax\ifmmode\else\unskip\hfill\fi\oprocendsymbol}

\begin{document}

\title{On the convergence time of asynchronous distributed quantized averaging algorithms}

\author{Minghui Zhu and Sonia Mart{\'\i}nez \thanks{The authors are
   with Department of Mechanical and Aerospace Engineering,
   University of California, San Diego, 9500 Gilman Dr, La Jolla CA,
   92093, {\tt\small \{mizhu,soniamd\}@ucsd.edu}}}

\maketitle

\begin{abstract} We come up with a class of distributed quantized averaging algorithms on asynchronous communication networks with fixed, switching and random topologies. The implementation of these algorithms is subject to the realistic constraint that the communication rate, the memory capacities of agents and the computation precision are finite. The focus of this paper is on the study of the convergence time of the proposed quantized averaging algorithms. By appealing to random walks on graphs, we derive polynomial bounds on the expected convergence time of the algorithms presented.
\end{abstract}

\section{Introduction} \label{sec:Introduction}

Consider a network of (mobile or immobile) agents. The distributed consensus problem aims to design an algorithm that agents can utilize to asymptotically reach an agreement by communicating with nearest neighbors. This problem historically roots in parallel computation~\cite{Bertsekas}, and has attracted significant attention recently~\cite{Blondel}\cite{Jadbabaie}\cite{Olfati1}. As a special case of the consensus problem, the distributed averaging problem requires that the consensus value be the average of individual
initial states. The distributed averaging algorithm acts as the building block for many distributed tasks such as sensor fusion~\cite{Xiao} and distributed estimation~\cite{Martinez}.


In real-world communication networks, the capacities of communication
channels and the memory capacities of agents are finite.
Furthermore, the computations can only be carried out with finite
precision. From a practical point of view, real-valued averaging
algorithms are not feasible and these realistic constraints motivate
the problem of average consensus via quantized information. Another
motivation for distributed quantized averaging is load balancing
with indivisible tasks. Prior work on distributed
quantized averaging over fixed graphs includes~\cite{Aysal,Carli1,Carli2,Kashyap}.
Recently,~\cite{Nedic} examines quantization effects on distributed averaging algorithms
over time-varying topologies. As in~\cite{Kashyap}, we focus on
quantized averaging algorithms preserving the sum of the state
values at each iteration. This setup has the following properties of
interest: the sum cannot be changed in some situations, such as
load balancing; and the constant sum leads to a small steady-state error with respect to the average of
individual initial states. This error is equal to either one
quantization step size or zero (when the average of the initial
states is located at one of the quantization levels) and thus is independent of
$N$. This is in contrast to the setup in~\cite{Nedic} where the sum
of the states is not maintained, resulting in a steady-state error of
the order $O(N^3\log N)$.


The convergence time is typically utilized to quantify the performance
of distributed averaging algorithms. The authors in~\cite{Boyd,Olshevsky}
study the convergence time of real-valued averaging, while~\cite{Kashyap,Nedic} discuss the
case of quantized averaging. The polynomial bounds of the expected convergence time on fixed complete
and linear graphs are derived in~\cite{Kashyap}. Recently, the authors in~\cite{Nedic}
give a polynomial bound on the convergence time of a class of quantized averaging algorithms
over switching topologies. Among the papers aforementioned, \cite{Boyd,Nedic,Olshevsky} require global synchronization, and \cite{Kashyap} needs some global information (e.g, a centralized entity or the total number of the edges) to explicitly bound the expected convergence times. However, real-world communication networks are inherently asynchronous environment and lack of centralized coordination.

\emph{Statement of contributions.} The present paper proposes a class of distributed
quantized averaging algorithms on asynchronous communication networks with fixed,
switching and random topologies. The algorithms are shown to asymptotically reach quantized average consensus in probability. Furthermore, we utilize meeting times of two random walks on
graphs as a unified approach to derive polynomial bounds on the expected convergence times of our presented algorithms. To the best of our knowledge, this note is the first step toward characterizing the expected convergence times of completely distributed quantized averaging algorithms over asynchronous communication networks. A preliminary conference version of this paper is in~\cite{Zhu} where the convergence time of synchronous algorithms is also studied.



\section{Preliminaries and problem statement} \label{sec:Problem}

Here, we present the problem formulation along with some notation and terminology.

\emph{Asynchronous time model.} In this note, we will employ the
asynchronous time model proposed in~\cite{Boyd}. More precisely, consider a network of $N$ nodes,
labeled $1$ through $N$. Each node has a clock which ticks
according to a rate $1$ Poisson process. Hence, the inter-tick
times at each node are random variables with rate $1$ exponential
distribution, independent across nodes and independent over time.
By the superposition theorem for Poisson processes, this setup is
equivalent to a single global clock modeled as a rate $N$ Poisson
process ticking at times $\{Z_k\}_{k\geq 0}$. By the orderliness
property of Poisson processes, the clock ticks do not occur
simultaneously. The inter-agent communication and the update of
consensus states only occur at $\{Z_k\}_{k\geq 0}$. In the reminder of this paper,
the time instant $t$ will be discretized according to $\{Z_k\}_{k\geq 0}$ and
defined in terms of the number of clock ticks.

\emph{Network model.} We will employ the
undirected graph ${\mathcal{G}}(t)=(V, E(t))$ to model the network. Here $V := \{1,\cdots,N\}$ is the vertex set, and an edge $(j,i)\in E(t)$ if and only if node $j$ can receive the message from node $i$
(e.g., node $j$ is within the communication range of node $i$) at time $t$. The neighbors of node $i$ at time $t$ are denoted by ${\mathcal{N}}_i(t)=\{j\in V \;|\; (j,i)\in E(t) \text{ and } j\neq i\}$. The state of node $i$ at time $t$ is denoted by $x_i(t)\in \mathbb{R}$ and the network state is denoted by
$x(t)=(x_1(t),\cdots,x_N(t))^T$. Suppose the initial states $x_i(0)\in[U_{\min}, U_{\max}]$
for all $i\in V$ and some real numbers $U_{\min}$ and $U_{\max}$.

\emph{Quantization scheme.} Let $R$ denote the number of bits per
sample. The total number of quantization levels can be represented by
$L=2^R$ and the step size is $\Delta=(U_{\max}-U_{\min})/2^R$. The set
of quantization levels, $\{\omega_1,\cdots,\omega_{L}\}$, is a
strictly increasing sequence in $\real$ and the levels are uniformly
spaced in the sense that $\omega_{i+1}-\omega_i=\Delta$. A quantizer
${\mathcal{Q}} : [U_{\min}, U_{\max}] \rightarrow
\{\omega_1,\cdots,\omega_{L}\}$ is adopted to quantize the message
$u\in[U_{\min}, U_{\max}]$ in such a way that
${\mathcal{Q}}(u)=\omega_i$ if $u\in[\omega_i,\omega_{i+1})$ for some
$i\in\{1,\cdots,L-1\}$. Assume that the initial states $x_i(0)$ for
all $i\in V$ are multiples of $\Delta$.

\emph{Problem statement.} The problem of interest in this paper is to
design distributed averaging algorithms which the nodes can utilize to update their
states by communicating with neighbors via quantized messages in an asynchronous setting. Ultimately, quantized average consensus is reached in probability; i.e., for any initial state $x(0)$, there holds that
$\lim_{t\rightarrow\infty}{\mathbb{P}}(x(t)\in {\mathcal{W}}(x(0)))=1$.  The
set ${\mathcal{W}}(x(0))$ is dependent on initial state $x(0)\in{\mathbb{R}}^N$ and defined as
follows. If $\bar{x}(0)=\frac{1}{N}\sum_{i=1}^N x_i(0)$ is not a
multiple of $\Delta$, then ${\mathcal{W}}(x(0)) = \{x\in {\mathbb{R}}^N\; |\;
x_i\in\{{\mathcal{Q}}(\bar{x}(0)), {\mathcal{Q}}(\bar{x}(0))+\Delta\}\}$;
otherwise, ${\mathcal{W}}(x(0)) = \{x\in {\mathbb{R}}^N\; |\; x_i=\bar{x}(0)\}$.
Now it is clear that the steady-state error is at most $\Delta$
after quantized average consensus is reached.

\emph{Notions of random walks on graphs.} In this paper, random walks on graphs play an
important role in characterizing the convergence properties of our quantized averaging algorithms.
The following definitions are generalized from those defined for fixed graphs
in~\cite{Brightwell,Coppersmith}.

\begin{definition}[Random walks] A random walk on the graph ${\mathcal{G}}(t)$ under the transition matrix
  $P(t)=(p_{ij}(t))$, starting from node $v$ at time $s$, is a stochastic process $\{X(t)\}_{t\geq s}$
  such that $X(s)=v$ and ${\mathbb{P}}(X(t+1)=j\;|\; X(t)=i)=p_{ij}(t)$. A random walk is said to be simple
  if for any $i\in V$, $p_{ii}(t)=0$ for all $t\geq0$; otherwise, it is said to be natural.\label{def1}\oprocend
\end{definition}


\begin{definition}[Hitting time] Consider a random walk on the graph
  ${\mathcal{G}}(t)$, beginning from node $i$ at time $s$ and evolving
  under the transition matrix $P(t)$. The hitting time from node $i$
  to the set $\Lambda\subseteq V$, denoted as
  $H_{({\mathcal{G}}(t),P(t),s)}(i,\Lambda)$, is the expected time it
  takes this random walk to reach the set $\Lambda$ for the first time.
  We denote
  $H_{({\mathcal{G}}(t),P(t))}(\Lambda)=\sup_{s\geq0}\max_{i\in V}H_{({\mathcal{G}}(t),P(t),s)}(i,\Lambda)$ as the hitting time to reach the set $\Lambda$. The hitting time of the pair $i,j$,
  denoted as $H_{({\mathcal{G}}(t),P(t),s)}(i,j)$, is the expected
  time it takes this random walk to reach node $j$ for the first time.
  Denote $H_{({\mathcal{G}}(t),P(t))}=\sup_{s\geq0}\max_{i,j\in V}H_{({\mathcal{G}}(t),P(t),s)}(i,j)$
  as the hitting time of going between any pair of nodes. \oprocend
\label{def2}\end{definition}

\begin{definition}[Meeting time] Consider two random walks on the
  graph ${\mathcal{G}}(t)$ under the transition matrix $P(t)$,
  starting at time $s$ from node $i$ and node $j$ respectively. The
  meeting time $M_{({\mathcal{G}}(t),P(t),s)}(i,j)$ of these two
  random walks is the expected time it takes them to meet at some node
  for the first time. The meeting time on the graph ${\mathcal{G}}(t)$
  is defined as $M_{({\mathcal{G}}(t),P(t))}=\sup_{s\geq0}\max_{i,j\in V}M_{({\mathcal{G}}(t),P(t),s)}(i,j)$.\oprocend
  \label{def4}\end{definition}

For the ease of notation, we will drop the subscript $s$ in the
hitting time and meeting time notions for fixed graphs. The
following notion is only defined for fixed graphs.

\begin{definition} [Irreducibility and reversibility]
A random walk on the graph $\mathcal{G}$ is irreducible if it is possible to get to any other node
from any node. An irreducible random walk with stationary distribution $\pi$ is called reversible if $\pi_i p_{ij}=\pi_j p_{ji}$ for all $i,j\in V$.\label{def5}\oprocend
\end{definition}

\emph{Notations.} For $\alpha\in {\mathbb{R}}$, define $V_{\alpha} : {\mathbb{R}}^N
\rightarrow {\mathbb{R}}$ as
$V_{\alpha}(x)=\sum_{i=1}^N(x_i-\alpha)^2$. We define $J : {\mathbb{R}}^N
\rightarrow {\mathbb{R}}$ as $J(x)=(\max_{i\in
V}x_i-\min_{i\in V}x_i)/\Delta$. Denote the set $\Theta=\{ (k,k) \;|\; k\in V\}$. The $distribution$ of a vector $x\in {\mathbb{R}}^N$ is defined to be the list $\{(q_1,m_1),(q_2,m_2),\cdots,(q_k,m_k)\}$ for some $k\in V$ where $\sum_{\ell=1}^k m_{\ell}=N$, $q_i\neq q_j$ for
$i\neq j$ and $m_{\ell}$ is the cardinality of the set $\{i\in V\;
|\;x_i=q_{\ell}\}$. The cardinality of the set $M$ is denoted by $|M|$.

\section{Asynchronous distributed quantized averaging on fixed graphs} \label{sec:AF}

In this section, we propose and analyze an asynchronous distributed quantized
averaging algorithm on the fixed and connected graph $\mathcal{G}$. Main references
are~\cite{Kashyap} on quantized gossip algorithms and ~\cite{Brightwell} on the meeting time of two simple random walks on fixed graphs.

\subsection{Proposed algorithm}

The \emph{asynchronous distributed quantized averaging algorithm on the fixed and
connected graph} $\mathcal{G}$ ($\AF$, for short) is described as
follows.  Suppose node $i$'s clock ticks at time $t$.
Node $i$ randomly chooses one of its neighbors, say node $j$, with
equal probability. Node $i$ and $j$ then
execute the following local computation. If $x_i(t)\geq x_j(t)$, then
\begin{align} x_i(t+1)=x_i(t)-\delta,\quad
x_j(t+1)=x_j(t)+\delta; \label{e156}\end{align} otherwise,
\begin{align} x_i(t+1)=x_i(t)+\delta,\quad
x_j(t+1)=x_j(t)-\delta,\label{e157}\end{align} where
$\delta=\frac{1}{2}(x_i(t)-x_j(t))$ if
$\frac{x_i(t)-x_j(t)}{2\Delta}$ is an integer; otherwise,
$\delta={\mathcal{Q}}(\frac{1}{2}(x_i(t)-x_j(t)))+\Delta$.
Every other node $k\in V\setminus \{i,j\}$ preserves its current state; i.e.,
$x_k(t+1)=x_k(t)$.

\begin{remark} The precision $\frac{\Delta}{2}$ is sufficient for the
  computation of $\delta$ and thus the update laws (\ref{e156}) and
  (\ref{e157}). It is easy to verify that $x_i(t)\in[U_{\min},
  U_{\max}]$ and $x_i(t)$ are multiples of $\Delta$ for all $i\in V$
  and $t\geq0$. Furthermore, the sum of the state values is preserved
  at each iteration.

  If $|x_i(t)-x_j(t)| = \Delta$, the update laws
  (\ref{e156}) and (\ref{e157}) become $x_i(t+1)=x_j(t)$ and
  $x_j(t+1)=x_i(t)$. Such update is referred to as a \emph{trivial
    average} in~\cite{Kashyap}. If $|x_i(t)-x_j(t)|>\Delta$, then
  (\ref{e156}) or (\ref{e157}) is referred to as a \emph{non-trivial
    average}.  Although it does not directly contribute to reaching
  quantized average consensus, trivial average helps the information
  flow over the network.\label{rem2}\oprocend\end{remark}

\subsection{The meeting time of two natural random walks on the fixed graph $\mathcal{G}$}\label{sec:meetingtimefixed}

To analyze the convergence properties of $\AF$, we first study a variation of the problem in~\cite{Coppersmith},
namely, \emph{the meeting time of two natural random walks on the fixed graph
$\mathcal{G}$}. More precisely, assume that the fixed graph $\mathcal{G}$ be undirected and connected.
Initially, two tokens are placed on the graph $\mathcal{G}$; at each time, one of the tokens is chosen with
probability $\frac{1}{N}$ and the chosen token moves to one of
the neighboring nodes with equal probability. What is the meeting
time for these two tokens?

The tokens move as two natural random walks with the transition matrix
$P_{\AF}$ on the graph $\mathcal{G}$. The matrix
$P_{\AF}=(\tilde{p}_{ij})\in {\mathbb{R}}^{N\times N}$ is given by
$\tilde{p}_{ii}=1-\frac{1}{N}$ for $i\in V$,
$\tilde{p}_{ij}=\frac{1}{N|{\mathcal{N}}_i|}$ for $(i,j)\in E$.
The meeting time of these two natural random walks is denoted as $M_{({\mathcal{G}},P_{\AF})}$.
Denote any of these two natural random walks as $X_{\mathcal{N}}$.
Correspondingly, we construct a simple random walk, say
$X_{\mathcal{S}}$, with the transition matrix $P_{\SF}$ on the graph
$\mathcal{G}$ where the matrix $P_{\SF}=(p_{ij})\in
{\mathbb{R}}^{N\times N}$ is given by $p_{ii}=0$ and
$p_{ij}=\frac{1}{|{\mathcal{N}}_i|}$ if $(i,j)\in E$. The
hitting times of the random walks $X_{\mathcal{S}}$ and $X_{\mathcal{N}}$ are denoted as
$H_{({\mathcal{G}},P_{\SF})}$ and $H_{({\mathcal{G}},P_{\AF})}$, respectively.

\begin{proposition} The meeting time of two natural random
  walks with transition matrices $P_{\AF}$ on the fixed graph
  $\mathcal{G}$ satisfies that $M_{({\mathcal{G}},P_{\AF})}\leq
  2NH_{({\mathcal{G}},P_{\SF})}-N$.
\label{pro1}\end{proposition}

\begin{proof}
Since the fixed graph $\mathcal{G}$ is undirected and connected, the random walks $X_{\mathcal{N}}$ and
$X_{\mathcal{S}}$ are irreducible. The reminder of the proof is based on the following claims:

(i) It holds that $H_{({\mathcal{G}},P_{\AF})}\geq N$.

(ii) For any pair $i,j\in V$ with $i\neq j$, we have $H_{({\mathcal{G}},P_{\AF})}(i,j)=NH_{({\mathcal{G}},P_{\SF})}(i,j)$.

(iii) For any $i,j,k\in V$, the following
equality holds:\\
$H_{({\mathcal{G}},P_{\AF})}(i,j)+H_{({\mathcal{G}},P_{\AF})}(j,k)
+H_{({\mathcal{G}},P_{\AF})}(k,i)=H_{({\mathcal{G}},P_{\AF})}(i,k)+H_{({\mathcal{G}},P_{\AF})}(k,j)
+H_{({\mathcal{G}},P_{\AF})}(j,i)$.

(iv) There holds that
$M_{({\mathcal{G}},P_{\AF})}\leq 2H_{({\mathcal{G}},P_{\AF})}-N.$

Now, let us prove each of the above claims.

(i) The quantity $H_{({\mathcal{G}},P_{\AF})}(i,j)$ reaches the
minimum when ${\mathcal{N}}_i=\{j\}$. We now consider the graph ${\mathcal{G}}$ with
${\mathcal{N}}_i=\{j\}$ and compute $H_{({\mathcal{G}},P_{\AF})}(i,j)$. The probability that
$X_{\mathcal{N}}$ stays up with node $i$ before time $\ell$ and
moves to node $j$ at time $\ell$ is $\frac{1}{N}(1-\frac{1}{N})^{\ell-1}$. Then, we have
$H_{({\mathcal{G}},P_{\AF})}(i,j)=\sum_{\ell=1}^{+\infty}\ell\frac{1}{N}(1-\frac{1}{N})^{\ell-1}=N$ and
Claim (i) holds.

(ii) For any pair $i,j\in V$ with $i\neq j$, it holds that
$H_{({\mathcal{G}},P_{\AF})}(i,j)=\sum_{k\in
{\mathcal{N}}_i}\frac{1}{N|\mathcal{N}_i|}(H_{({\mathcal{G}},P_{\AF})}(k,j)+1)
+(1-\frac{1}{N})(H_{({\mathcal{G}},P_{\AF})}(i,j)+1)$. Hence, we have
that $H_{({\mathcal{G}},P_{\AF})}(i,j)=
N+\sum_{k\in{{\mathcal{N}}_i}}\frac{1}{|\mathcal{N}_i|}H_{({\mathcal{G}},P_{\AF})}(k,j)$.
Furthermore,
$H_{({\mathcal{G}},P_{\SF})}(i,j)=\sum_{k\in
{\mathcal{N}}_i}\frac{1}{|\mathcal{N}_i|}(H_{({\mathcal{G}},P_{\SF})}(k,j)+1)=1+\sum_{k\in
{\mathcal{N}}_i}\frac{1}{|\mathcal{N}_i|}H_{({\mathcal{G}},P_{\SF})}(k,j)$.
Hence, Claim (ii) holds.

(iii) Denote by $\pi_i=|{\mathcal{N}}_i|/{\mathcal{N}}_{\max}$ and
$\pi=(\pi_1,\cdots,\pi_N)^T$ where ${\mathcal{N}}_{\max}=\max_{i\in
V}\{|{\mathcal{N}}_i|\}$. Since $P_{\AF}^T\pi=\pi$, then $\pi$ is
the stationary distribution of the random walk $X_{\mathcal{N}}$. Furthermore, for
any pair $i,j\in V$, we have $\pi_i
\tilde{p}_{ij}=\frac{|{\mathcal{N}}_i|}{{\mathcal{N}}_{\max}}\frac{1}{N|{\mathcal{N}}_i|}
=\frac{1}{N{\mathcal{N}}_{\max}}=\pi_j
\tilde{p}_{ji}=\frac{|{\mathcal{N}}_j|}{{\mathcal{N}}_{\max}}\frac{1}{N|{\mathcal{N}}_j|}$
and thus the random walk $X_{\mathcal{N}}$ is reversible. From Lemma 2 of
\cite{Coppersmith} it follows that Claim (iii) holds.

(iv) Claim (iv) is an extension of Theorem 2 in \cite{Coppersmith}.
An immediate result of Claim (iii) gives a node-relation on $V$;
i.e., $i\leq j$ if and only if $H_{({\mathcal{G}},P_{\AF})}(i,j)\leq
H_{({\mathcal{G}},P_{\AF})}(j,i)$. This relation is transitive and
constitutes a pre-order on $V$. Then there exists a node $u$
satisfying $H_{({\mathcal{G}},P_{\AF})}(v,u)\geq
H_{({\mathcal{G}},P_{\AF})}(u,v)$ for any other node $v\in V$. Such
a node $u$ is called $hidden$. As in~\cite{Coppersmith}, we define a
potential function $\Phi$ by
$\Phi(i,j)=H_{({\mathcal{G}},P_{\AF})}(i,j)+H_{({\mathcal{G}},P_{\AF})}(j,u)-H_{({\mathcal{G}},P_{\AF})}(u,j)$.

Define the functions $\Phi(\bar{i},j)$ and
$M_{({\mathcal{G}},P_{\AF})}(\bar{i},j)$ below, the averages of the
functions $\Phi$ and $M_{({\mathcal{G}},P_{\AF})}$ over the
neighbors of node $i$ and $j$, respectively:
\begin{align*}&\Phi(\bar{i},j)=\frac{1}{|\mathcal{N}_i|}\sum_{k\in
{{\mathcal{N}}_i}}\Phi(k,j)=\frac{1}{|\mathcal{N}_i|}\sum_{k\in
{{\mathcal{N}}_i}}H_{({\mathcal{G}},P_{\AF})}(k,j)
+H_{({\mathcal{G}},P_{\AF})}(j,u)-H_{({\mathcal{G}},P_{\AF})}(u,j),\nnum\\
&M_{({\mathcal{G}},P_{\AF})}(\bar{i},j)=\frac{1}{|\mathcal{N}_i|}
\sum_{k\in{{\mathcal{N}}_i}}M_{({\mathcal{G}},P_{\AF})}(k,j).\end{align*}

In Claim (ii), we have shown that
$H_{({\mathcal{G}},P_{\AF})}(i,j)=\sum_{k\in{{\mathcal{N}}_i}}\frac{1}{|\mathcal{N}_i|}
H_{({\mathcal{G}},P_{\AF})}(k,j)+N$. Thus, $\Phi(\bar{i},j)+N=\Phi(i,j)$. Similarly,
$M_{({\mathcal{G}},P_{\AF})}(\bar{i},j)+N=M_{({\mathcal{G}},P_{\AF})}(i,j)$.

We are now in a position to show that for any pair $i,j\in V$, it
holds that \begin{align}M_{({\mathcal{G}},P_{\AF})}(i,j)\leq
\Phi(i,j).\label{e197}\end{align} Assume that (\ref{e197}) does not
hold. Let $\phi$ be $\phi=\max_{w,v\in
V}(M_{({\mathcal{G}},P_{\AF})}(w,v)-\Phi(w,v))>0$. Choose a pair of
$i,j$ with minimum distance among the set $\Xi=\{(w,v)\in V\times
V\;|\; M_{({\mathcal{G}},P_{\AF})}(w,v)-\Phi(w,v)=\phi\}$. Toward
this end, consider the following two cases:

(1) $j\in {\mathcal{N}}_i$. Observe that
$\Phi(j,j)=H_{({\mathcal{G}},P_{\AF})}(j,j)+H_{({\mathcal{G}},P_{\AF})}(j,u)
-H_{({\mathcal{G}},P_{\AF})}(u,j)\geq0
=M_{({\mathcal{G}},P_{\AF})}(j,j)$. We have
$\Phi(\bar{i},j)+\phi>M_{({\mathcal{G}},P_{\AF})}(\bar{i},j)$ and
thus
\begin{align}
M_{({\mathcal{G}},P_{\AF})}(i,j)=\Phi(i,j)+\phi=N+\Phi(\bar{i},j)+\phi>N+M_{({\mathcal{G}},P_{\AF})}(\bar{i},j)
=M_{({\mathcal{G}},P_{\AF})}(i,j).\label{e195}
\end{align}

(2) $j\notin {\mathcal{N}}_i$. There exists node $k\in
{\mathcal{N}}_i$ such that node $k$ is closer to node $j$ than node
$i$. Since the pair of $i,j$ has the minimum distance in the set
$\Xi$, we have $M_{({\mathcal{G}},P_{\AF})}(k,j)-\Phi(k,j)<\phi$.
Hence,
$\Phi(\bar{i},j)+\phi>M_{({\mathcal{G}},P_{\AF})}(\bar{i},j)$, and
thus \eqref{e195} holds.

In both cases, we get to the contradiction
$M_{({\mathcal{G}},P_{\AF})}(i,j)>M_{({\mathcal{G}},P_{\AF})}(i,j)$, and
thus \eqref{e197} holds.

Combining Claims (i), (ii) and inequality (\ref{e197}) gives the desired result of
$M_{({\mathcal{G}},P_{\AF})}\leq2NH_{({\mathcal{G}},P_{\SF})}-N$.
\end{proof}

\subsection{Convergence analysis of $\AF$}

We now proceed to analyze the convergence properties of $\AF$. The convergence
time of $\AF$ is a random variable defined as follows:
$T_{\rm con}(x(0))=\inf\{t\;|\; x(t)\in {\mathcal{W}}(x(0))\},$
where $x(t)$ starts from $x(0)$ and evolves under $\AF$. Choose $V_{\bar{x}(0)}(x)=\sum_{i=1}^N(x_i-\bar{x}(0))^2$ as a Lyapunov function candidate for $\AF$. One can readily see that
$V_{\bar{x}(0)}(x(t+1))=V_{\bar{x}(0)}(x(t))$ when a trivial average
occurs and $V_{\bar{x}(0)}(x)$ reduces at least $2\Delta^2$ when a
non-trivial average occurs. Hence, $V_{\bar{x}(0)}(x)$ is
non-increasing along the trajectories, and the number of non-trivial averages is at most
$\frac{1}{2\Delta^2}V_{\bar{x}(0)}(x(0))$. Define the set $\Psi=\{x \in
{\mathbb{R}}^N\;|\; {\rm the~ distribution~ of~} x {\rm ~is~}
\{(0,1),(\Delta,N-2),(2\Delta,1)\}\}$ and denote
${\mathbb{E}}[T_{\Psi}]=\max_{x(0)\in \Psi}{\mathbb{E}}[T_{\rm
con}(x(0))]$. It is clear that the expected time between any two consecutive
non-trivial averages is not larger than ${\mathbb{E}}[T_{\Psi}]$.
Then we have the following estimates on ${\mathbb{E}}[T_{\rm con}(x(0))]$:
\begin{align}{\mathbb{E}}[T_{\rm con}(x(0))]\leq\frac{1}{2\Delta^2}V_{\bar{x}(0)}(x(0)){\mathbb{E}}[T_{\Psi}]\leq
  \frac{N J(x(0))^2}{8}{\mathbb{E}}[T_{\Psi}], \label{e106}\end{align} where
the second inequality is a direct result of Lemma 4 in~\cite{Kashyap}.

\begin{theorem} For any $x(0)\notin {{\mathcal{W}}}(x(0))$, the expected convergence time
${\mathbb{E}}[T_{\rm con}(x(0))]$ of $\AF$ is upper bounded by $\frac{N^2
J(x(0))^2}{8}(\frac{8}{27}N^3-1)$. \label{the2}
\end{theorem}

\begin{proof} By (\ref{e106}), it suffices to bound ${\mathbb{E}}[T_{\Psi}]$. Assume that
  $x(0)\in \Psi$. Before they meet for the first time, the values
  $0$ and $2\Delta$ move as two natural random walks which are
  identical to $X_{\mathcal{N}}$ in Proposition~\ref{pro1}. At their meeting for the first time, the values of $0$ and $2\Delta$ average and quantized average consensus is reached. Hence, ${\mathbb{E}}[T_{\Psi}]=M_{({\mathcal{G}},P_{\AF})}$ and thus inequality~\eqref{e106} becomes
  \begin{align}
    {\mathbb{E}}[T_{\rm{con}}(x(0))]\leq \frac{N
      J(x(0))^2}{8}M_{({\mathcal{G}},P_{\AF})}\leq \frac{N
      J(x(0))^2}{8}(2NH_{({\mathcal{G}},P_{\SF})}-N),\label{e101}
  \end{align}
  where we use Proposition~\ref{pro1} in the second inequality. By letting $M=0$ in the theorem of Page 265 in~\cite{Brightwell}, we can obtain the upper bound $\frac{4}{27}N^3$ on $H_{({\mathcal{G}},P_{\SF})}$. Substituting this upper bound into inequality~\eqref{e101} gives the desired upper bound on ${\mathbb{E}}[T_{\rm con}(x(0))]$.
\end{proof}

\begin{theorem} Let $x(0)\in {\mathbb{R}}^N$ and suppose $x(0)\notin {\mathcal{W}}(x(0))$.
Under $\AF$, almost any evolution $x(t)$ starting from $x(0)$ reaches quantized average consensus.\label{the4}
\end{theorem}

\begin{proof} Denote $\tilde{T} = \frac{N^2J(x(0))^2}{4}(\frac{8}{27}N^3-1)$, and
consider the first $\tilde{T}$ clock ticks of evolution of $\AF$ starting from $x(0)$. It follows
from Markov's inequality that \begin{align*} {\mathbb{P}}(T_{\rm
con}(x(0))>\tilde{T}\;|\;x(0)\notin{\mathcal{W}}(x(0)))\leq\frac{{\mathbb{E}}[T_{\rm
con}(x(0))]}{\tilde{T}}\leq\frac{1}{2},\end{align*} that is, the probability that after $\tilde{T}$
clock ticks $\AF$ has not reached quantized average consensus is less than $\frac{1}{2}$.
Starting from $x(\tilde{T})$, let us consider the posterior
evolution of $x(t)$ in the next $\tilde{T}$ clock ticks. We have
\begin{align*}
{\mathbb{P}}(T_{\rm con}(x(\tilde{T}))>\tilde{T}\;|\;x(\tilde{T})\notin
{\mathcal{W}}(x(0)))\leq\frac{{\mathbb{E}}[T_{\rm
con}(x(\tilde{T}))]}{\tilde{T}} \leq\frac{1}{2}. \end{align*} That
is, the probability that after $2\tilde{T}$ clock ticks $x(t)$ has
not reached quantized average consensus is at most
$(\frac{1}{2})^2$. By induction, it follows that after $n\tilde{T}$
clock ticks the probability $x(t)$ not reaching quantized average
consensus is at most $(\frac{1}{2})^n$. Since the set ${\mathcal{W}}(x(0))$ is absorbing, we have $\lim_{t\rightarrow\infty}{\mathbb{P}}(x(t)\notin {\mathcal{W}}(x(0)))=0$.
This completes the proof.\end{proof}

\section{Asynchronous distributed quantized averaging on switching graphs}\label{sec:AS}

We now turn our attention to the more challenging scenario where the communication graphs
are undirected but dynamically changing. We will propose and analyze an
\emph{asynchronous distributed quantized averaging algorithm on switching graphs}
($\AS$, for short). The convergence rate of distributed real-valued averaging algorithms on switching graphs
in~\cite{Nedic} will be employed to characterize the hitting time of random walks on switching graphs.

\subsection{Proposed algorithm}

The main steps of $\AS$ can be summarized as follows. At time $t$, let node $i$'s clock tick. If
$|{\mathcal{N}}_i(t)|\neq0$, node $i$ randomly chooses one of its
neighbors, say node $j$, with probability
$\frac{1}{\max\{|{\mathcal{N}}_i(t)|,|{\mathcal{N}}_j(t)|\}}$. Then,
node $i$ and $j$ execute the computation (\ref{e156}) or
(\ref{e157}) and every other node $k\in V\backslash\{i,j\}$ preserves its current state.
If $|{\mathcal{N}}_i(t)|=0$, all the nodes do nothing at this time.

Here, we assume that the communication graph ${\mathcal{G}}(t)$ be undirected and satisfies the
following connectivity assumption also used in
~\cite{Blondel,Jadbabaie,Nedic,Olshevsky}.


\begin{assumption} [Periodical connectivity] There exists some $B\in{\mathbb{N}}_{>0}$ such that, for all $t\geq0$, the undirected graph $(V,E(t)\cup E(t+1)\cup\cdots\cup E(t+B-1))$ is connected.\label{asm12}\end{assumption}

\begin{remark} In the $\AS$, the
  probability that node $i$ chooses a neighbor $j$ is
  $\frac{1}{\max\{|{\mathcal{N}}_i(t)|,|{\mathcal{N}}_j(t)|\}}$. Thus,
  this information should be available to node $i$. In this way, the
  matrix $P_{\AS}(t)$ defined later is symmetric and double
  stochastic.\oprocend\label{rem3}\end{remark}

\subsection{The meeting time of two natural random walks on the time-varying graph ${\mathcal{G}}(t)$}

Before analyzing $\AS$, we consider the following problem which generalizes the problem in Section~\ref{sec:meetingtimefixed} to the case of dynamically changing graphs.

\emph{The meeting time of two natural random walks on the time-varying graph ${\mathcal{G}}(t)$}.
Assume that ${\mathcal{G}}(t)$ be undirected and satisfies Assumption~\ref{asm12}.
Initially, two tokens are placed on ${\mathcal{G}}(0)$. At each time, one of the tokens is chosen
with probability $\frac{1}{N}$. The chosen token at some node, say
$i$, moves to one of the neighbors, say node $j$, with probability
$\frac{1}{\max\{|{\mathcal{N}}_i(t)|,|{\mathcal{N}}_j(t)|\}}$ if
$|{\mathcal{N}}_i(t)|\neq0$; otherwise, it will stay up with node $i$. What is
the meeting time for these two tokens?

Clearly, the movements of two tokens are two natural random walks, say
$X_1$ and $X_2$, on the switching graph ${\mathcal{G}}(t)$. Their
meeting time is denoted as $M_{({\mathcal{G}}(t),P_{\AS}(t))}$ where
the transition matrix $P_{\AS}(t)=(\bar{p}_{ij}(t))$ is given as follows: if
$|{\mathcal{N}}_i(t)|\neq0$, then
$\bar{p}_{ij}(t)=\frac{1}{N\max\{|{\mathcal{N}}_i(t)|,|{\mathcal{N}}_j(t)|\}}$
for $(i,j)\in E(t)$ and $\bar{p}_{ii}(t)=1-\sum_{(i,j)\in
  E(t)}\frac{1}{N\max\{|{\mathcal{N}}_i(t)|,|{\mathcal{N}}_j(t)|\}}$;
if $|{\mathcal{N}}_i(t)|=0$, then $\bar{p}_{ii}(t)=1$. One can easily verify that the
matrix $P_{\AS}(t)$ is symmetric and doubly stochastic. The
natural random walks $X_1$ and $X_2$ on the graph ${\mathcal{G}}(t)$ are
equivalent to a single natural random walk, say $X_M$, on the product
graph ${\mathcal{G}}(t)\times {\mathcal{G}}(t)$. That is, $X_M$ moving
from node $(i_1,i_2)\in V \times V$ to node $(j_1,j_2)\in V \times V$ on the graph
${\mathcal{G}}(t)\times {\mathcal{G}}(t)$ at time $t$, is equivalent
to $X_1$ moving from $i_1$ to $j_1$ and $X_2$ moving from $i_2$ to
$j_2$ on the graph ${\mathcal{G}}(t)$ at time $t$. Denote the
transition matrix of the random walk $X_M$ as $Q(t)=(q_{(i_1,i_2)(j_1,j_2)}(t))\in
\mathbb{R}^{N^2\times N^2}$.

In the following lemma, we will consider the random walk $\bar{X}_M$ on the graph ${\mathcal{G}}(t)\times {\mathcal{G}}(t)$ with the absorbing set $\Theta$ and the transition matrix $\bar{Q}(t)\in
\mathbb{R}^{N^2\times N^2}$. Denote $e_{(\ell_1,\ell_2)}$ by the row corresponding to $(\ell_1,\ell_2)\in V\times V$ in a $N^2\times N^2$ identity matrix. The transition matrix $\bar{Q}(t)$ is defined by replacing the row associated with the absorbing state $(\ell_1,\ell_2)\in\Theta$ in $Q(t)$ with $e_{(\ell_1,\ell_2)}$. Define $\vartheta_{(\ell_1,\ell_2)}(t)=\mathbb{P}(X_M(t)=(\ell_1,\ell_2))$,
$\vartheta(t)=\col\{\vartheta_{(\ell_1,\ell_2)}(t)\}\in{\mathbb{R}}^{N^2}$,
$\vartheta_{\Theta}(t)=\sum_{(\ell_1,\ell_2)\in \Theta}\vartheta_{(\ell_1,\ell_2)}(t)$ for the random walk $X_M$, and
$\bar{\vartheta}_{(\ell_1,\ell_2)}(t)=\mathbb{P}(\bar{X}_M(t)=(\ell_1,\ell_2))$,
$\bar{\vartheta}(t)=\col\{\bar{\vartheta}_{(\ell_1,\ell_2)}(t)\}\in{\mathbb{R}}^{N^2}$,
$\bar{\vartheta}_{\Theta}(t)=\sum_{(\ell_1,\ell_2)\in \Theta}
\bar{\vartheta}_{(\ell_1,\ell_2)}(t)$ for the random walk $\bar{X}_M$.

\begin{lemma} Consider a network of $N$ nodes whose communication
  graph ${\mathcal{G}}(t)$ be undirected and satisfies Assumption~\ref{asm12}. Let $(i_1, i_2)\in
  V\times V$ be a given node and suppose that the random walks $X_M$ and $\bar{X}_M$ start from node
  $(i_1,i_2)$ at time $0$. Then it holds that
  $\bar{\vartheta}_{\Theta}(t)\geq \vartheta_{\Theta}(t)\geq \frac{1}{2N}$ for $t\geq t_1$ where $t_1$ is the smallest integer which is larger than $B(8N^6\log (\sqrt{2}N)+1)$.\label{lem4}
\end{lemma}

\begin{proof} It is not difficult to check that ${\mathcal{G}}(t)\times {\mathcal{G}}(t)$
is undirected and satisfies Assumption~\ref{asm12} with period $B$.
The minimum of nonzero entries in $Q(t)$ is lower bounded by $\frac{1}{N(N-1)}$, and
$Q(t)$ is symmetric. Observe that for any $(i_1,i_2)\in  V\times V$ and any $t\geq0$, $\sum_{(j_1,j_2)\in V\times V}
q_{(i_1,i_2)(j_1,j_2)}(t) = \sum_{(j_1,j_2)\in V\times V}
\bar{p}_{i_1 j_1}(t)\bar{p}_{i_2 j_2}(t) = \sum_{j_1\in V}
\bar{p}_{i_1 j_1}(t) \times \sum_{j_2\in V}
\bar{p}_{i_2 j_2}(t) = 1$ where we use the fact that the matrix $P_{\AS}(t)$ is doubly stochastic.
Hence, the matrix $Q(t)$ is doubly stochastic.

The evolution of $\vartheta(t)$ is governed by the
equation $\vartheta(t+1) = Q^T(t)\vartheta(t)$ with initial
state $\vartheta(0)=e_{(i_1,i_2)}^T$. Consider the Lyapunov function
$V_{\frac{1}{N^2}}(\vartheta)=\sum_{i=1}^{N^2}(\vartheta_i-\frac{1}{N^2})^2$
with $V_{\frac{1}{N^2}}(\vartheta(0))=1-\frac{1}{N^2}$.
It follows from Lemma 5 in~\cite{Nedic} that
\begin{align}
  V_{\frac{1}{N^2}}(\vartheta((k+1)B))\leq
  (1-\frac{1}{2N^5(N-1)})V_{\frac{1}{N^2}}(\vartheta(kB))\label{e112}
\end{align}
for $k\in {\mathbb{N}}_0$. Denote $\textbf{1}\in{\mathbb{R}}^{N^2}$ as the
vector of $N^2$ ones and note that
\begin{align*}
  V_{\frac{1}{N^2}}(\vartheta(t))-V_{\frac{1}{N^2}}(\vartheta(t+1))
  =(\vartheta(t)-\frac{1}{N^2}\textbf{1})^T(I-Q(t)Q^T(t))(\vartheta(t)-\frac{1}{N^2}\textbf{1}).
\end{align*}

Since $Q(t)$ is doubly stochastic, so is $Q(t)Q^T(t)$. Hence, the diagonal
entries of the matrix $\Gamma(t)=I-Q(t)Q^T(t)=(\gamma_{ij}(t))\in {\mathbb{R}}^{N^2\times N^2}$
are dominant in the sense of $\gamma_{ii}(t)=\sum_{j\neq i}\gamma_{ij}(t)$.  According
to Gershgorin theorem in~\cite{Horn}, all eigenvalues of $\Gamma(t)$
lie in a closed disk centered at $\max_{i\in \{1,\cdots,N^2\}}{\gamma_{ii}(t)}$ with
a radius $\max_{i\in \{1,\cdots,N^2\}}{\gamma_{ii}(t)}$.  Hence, $\Gamma(t)$ is
positive semi-definite.  Consequently,
$V_{\frac{1}{N^2}}(\vartheta(t))-V_{\frac{1}{N^2}}(\vartheta(t+1))\geq0$
and thus $V_{\frac{1}{N^2}}(\vartheta(t))$ is non-increasing along the trajectory of $\vartheta(t)$.
Combining (\ref{e112}) with the non-increasing property of
$V_{\frac{1}{N^2}}(\vartheta(t))$ gives that
\begin{align} V_{\frac{1}{N^2}}(\vartheta(t))\leq
V_{\frac{1}{N^2}}(\vartheta(0))(1-\frac{1}{2N^5(N-1)})^{\frac{t}{B}-1}
=\frac{N^2-1}{N^2}(1-\frac{1}{2N^5(N-1)})^{\frac{t}{B}-1}.\label{e114}
\end{align}

Since $\vartheta(t)^T \textbf{1} = 1$, then $\vartheta_{\min}(t) := \min_{(\ell_1,\ell_2)\in
V\times V}\vartheta_{(\ell_1,\ell_2)}(t)\leq \frac{1}{N^2}$. Since $V_{\frac{1}{N^2}}(\vartheta(t))\geq (\vartheta_{\min}(t)-\frac{1}{N^2})^2$, inequality (\ref{e114}) gives
that
$\vartheta_{\min}(t)\geq\frac{1}{N^2}-
(\frac{N^2-1}{N^2}(1-\frac{1}{2N^5(N-1)})^{\frac{t}{B}-1})^{\frac{1}{2}}$.
 Therefore, it holds that $\vartheta_{\min}(t)\geq
\frac{1}{2N^2}$ for $t\geq
B(\frac{\log(4N^2(N^2-1))}{-\log(1-\frac{1}{2N^5(N-1)})}+1)$. Since
$\log x\leq x-1$, there holds
$\frac{1}{-\log(1-\frac{1}{2N^5(N-1)})}\leq 2N^5(N-1)\leq 2N^6$. Hence, we
have that $\vartheta_{\min}(t)\geq \frac{1}{2N^2}$ and thus $\vartheta_{\Theta}(t)\geq \frac{1}{2N}$ for $t\geq t_1$.

Note that the evolution of $\bar{\vartheta}(t)$ is governed by the
equation $\bar{\vartheta}(t+1) = \bar{Q}(t)^T\bar{\vartheta}(t)$
with $\bar{\vartheta}(0)=e_{(i_1,i_2)}$. Since the set $\Theta$ is absorbing, $\bar{\vartheta}_{\Theta}(t)\geq \vartheta_{\Theta}(t)$ for all $t\geq0$ and thus the desired result follows.
\end{proof}

\begin{proposition} The meeting time of two
natural random walks with transition matrix $P_{\AS}(t)$ on the time-varying graph ${\mathcal{G}}(t)$ satisfies that
$M_{({\mathcal{G}}(t),P_{\AS}(t))}\leq 4Nt_1$.\label{pro2}
\end{proposition}

\begin{proof} Denote by $H_{({\mathcal{G}}(t)\times{\mathcal{G}}(t),Q(t))}(\Theta)$ the hitting time of
 the random walk $X_M$ to reach the set of $\Theta$. Observe that $M_{({\mathcal{G}}(t),P_{\AS}(t))}=H_{({\mathcal{G}}(t)\times{\mathcal{G}}(t),Q(t))}(\Theta)$. To find an upper bound on $H_{({\mathcal{G}}(t)\times{\mathcal{G}}(t),Q(t))}(\Theta)$, we construct the random
walk $X_M^{(i_1,i_2)}$ in such a way that $X_M^{(i_1,i_2)}$ starts from $(i_1,i_2)$
at time $0$ with $i_1 \neq i_2$ and the set $\Theta$ is the absorbing set of $X_M^{(i_1,i_2)}$. The transition matrix of $X_M^{(i_1,i_2)}$ is $\bar{Q}(t)$ defined before Lemma~\ref{lem4}. Define
$\vartheta_{(\ell_1,\ell_2)}^{(i_1,i_2)}(t)={\mathbb{P}}(X_M^{(i_1,i_2)}(t) = (\ell_1,\ell_2))$, and
$\vartheta^{(i_1,i_2)}(t)=\col\{\vartheta_{(\ell_1,\ell_2)}^{(i_1,i_2)}(t)\}\in
{\mathbb{R}}^{N^2}$. The dynamics of $\vartheta^{(i_1,i_2)}(t)$ is given by
$\vartheta^{(i_1,i_2)}(t+1)=\bar{Q}(t)^T\vartheta^{(i_1,i_2)}(t)$ with
the initial state $\vartheta^{(i_1,i_2)}(0)=e_{(i_1,i_2)}^T$.

Define the function 
$\mu_{(\ell_1,\ell_2)}^{(i_1,i_2)}: {\mathbb{N}}_{0}\rightarrow \{0,1\}$ in such a
way that $\mu_{(\ell_1,\ell_2)}^{(i_1,i_2)}=1$ if $X_M^{(i_1,i_2)}(t)=(\ell_1,\ell_2)$; otherwise,
$\mu_{(\ell_1,\ell_2)}^{(i_1,i_2)}(t)=0$. Define
$n_{(\ell_1,\ell_2)}^{(i_1,i_2)}=\sum_{\tau=0}^{+\infty}\mu_{(\ell_1,\ell_2)}^{(i_1,i_2)}(\tau)$ which
is the total times that the random walk $X_M^{(i_1,i_2)}$ is at node $(\ell_1,\ell_2)$. Then, the
hitting time $H_{({\mathcal{G}}(t)\times{\mathcal{G}}(t),Q(t),0)}((i_1,i_2),\Theta)$ of
$X_M^{(i_1,i_2)}$ equals the expected time that $X_M^{(i_1,i_2)}$ stays up with the nodes in $V\times V\backslash \Theta$, that is,
\begin{align}
  &H_{({\mathcal{G}}(t)\times{\mathcal{G}}(t),Q(t),0)}((i_1,i_2),\Theta) = \sum_{(\ell_1,\ell_2)\notin \Theta}{\mathbb{E}}[n_{(\ell_1,\ell_2)}^{(i_1,i_2)}]
  =\sum_{(\ell_1,\ell_2)\notin\Theta}{\mathbb{E}}[\sum_{\tau=0}^{+\infty}\mu_{(\ell_1,\ell_2)}^{(i_1,i_2)}(\tau)]\nnum\\
  &=\sum_{(\ell_1,\ell_2)\notin \Theta}\sum_{\tau=0}^{+\infty}{\mathbb{E}}[\mu_{(\ell_1,\ell_2)}^{(i_1,i_2)}(\tau)]
  =\sum_{\tau=0}^{+\infty}\sum_{(\ell_1,\ell_2)\notin \Theta}
  \vartheta_{(\ell_1,\ell_2)}^{(i_1,i_2)}(\tau).\label{e178}
\end{align}

It follows from Lemma~\ref{lem4} that $\vartheta_{\Theta}^{(i_1,i_2)}(t)\geq
\frac{1}{2N}$ for $t\geq t_1$. With that, the fact of $\vartheta^{(i_1,i_2)}(t)^T\textbf{1}=1$ implies that
\begin{align}\sum_{(\ell_1,\ell_2)\notin \Theta}
  \vartheta_{(\ell_1,\ell_2)}^{(i_1,i_2)}(t_1)\leq1-\frac{1}{2N}.\label{e177}\end{align}

For each $(k_1,k_2)\notin \Theta$, we construct the random walk $\tilde{X}_M^{(k_1,k_2)}$ in
such a way that $\tilde{X}_M^{(k_1,k_2)}$ starts from $(k_1,k_2)$ at time
$t_1$ and the set $\Theta$ is the absorbing set of
$\tilde{X}_M^{(k_1,k_2)}$. The transition matrix of $\tilde{X}_M^{(k_1,k_2)}$ is
$\bar{Q}(t)$. Define $\tilde{\vartheta}_{(\ell_1,\ell_2)}^{(k_1,k_2)}(t)=\mathbb{P}(\tilde{X}_M^{(k_1,k_2)}(t)=(\ell_1,\ell_2))$.
Following the forgoing arguments for $X^{(i_1,i_2)}_M$, we have
\begin{align}\sum_{(\ell_1,\ell_2)\notin \Theta}
  \tilde{\vartheta}_{(\ell_1,\ell_2)}^{(k_1,k_2)}(2t_1)\leq1-\frac{1}{2N}.\label{e179}\end{align}

Combining (\ref{e177}) and (\ref{e179}) gives that
\begin{align}&\sum_{(\ell_1,\ell_2)\notin \Theta}\vartheta_{(\ell_1,\ell_2)}^{(i_1,i_2)}(2t_1)
=\sum_{(\ell_1,\ell_2)\notin \Theta}\sum_{(k_1,k_2)\notin \Theta}\vartheta^{(i_1,i_2)}_{(k_1,k_2)}(t_1)\tilde{\vartheta}_{(\ell_1,\ell_2)}^{(k_1,k_2)}(2t_1)\nnum\\
&=\sum_{(k_1,k_2)\notin \Theta}\vartheta^{(i_1,i_2)}_{(k_1,k_2)}(t_1)
\sum_{(\ell_1,\ell_2)\notin \Theta}\tilde{\vartheta}^{(k_1,k_2)}_{(\ell_1,\ell_2)}(2t_1)\leq
(1-\frac{1}{2N})^2.\label{e184}\end{align} By induction, we have $\sum_{(\ell_1,\ell_2) \notin \Theta}\vartheta_{(\ell_1,\ell_2)}^{(i_1,i_2)}(nt_1)\leq (1-\frac{1}{2N})^n$ and then obtain
a strictly decreasing sequence $\sum_{(\ell_1,\ell_2) \notin \Theta}\vartheta_{(\ell_1,\ell_2)}^{(i_1,i_2)}(nt_1)$ with respect to $n\in{\mathbb{Z}}_0$. Since the set $\Theta$ is absorbing, then $\sum_{(\ell_1,\ell_2) \notin \Theta}\vartheta_{(\ell_1,\ell_2)}^{(i_1,i_2)}(t)$ is non-increasing with respect to $t\geq0$. Therefore, we have the following estimate
\begin{align}\sum_{(\ell_1,\ell_2) \notin \Theta}\vartheta_{(\ell_1,\ell_2)}^{(i_1,i_2)}(t)\leq
\sum_{(\ell_1,\ell_2) \notin \Theta}\vartheta_{(\ell_1,\ell_2)}^{(i_1,i_2)}(0)(1-\frac{1}{2N})^{\frac{t}{t_1}-1}
=(1-\frac{1}{2N})^{\frac{t}{t_1}-1}.\label{e180}\end{align}
Substituting (\ref{e180}) into (\ref{e178}) gives that
\begin{align}H_{({\mathcal{G}}(t)\times {\mathcal{G}}(t),Q(t),0)}((i_1,i_2),\Theta)
\leq \sum_{\tau=0}^{+\infty}
(1-\frac{1}{2N})^{\frac{\tau}{t_1}-1}=(1-\frac{1}{2N})^{-\frac{1}{t_1}}\cdot
\frac{1}{1-(1-\frac{1}{2N})^{\frac{1}{t_1}}}.\label{e102}\end{align}

Since $t_1>1$, it holds that $(1-\frac{1}{2N})^{-\frac{1}{t_1}}\leq
2^{\frac{1}{t_1}}<2$. It follows from Bernoulli's inequality that
$(1-\frac{1}{2N})^{\frac{1}{t_1}}\leq 1-\frac{1}{2Nt_1}$, and thus
$\frac{1}{1-(1-\frac{1}{2N})^{\frac{1}{t_1}}}\leq 2N t_1$.
Inequality (\ref{e102}) becomes
\begin{align}H_{({\mathcal{G}}(t)\times {\mathcal{G}}(t),Q(t),0)}((i_1,i_2),\Theta) \leq
4Nt_1.\label{e107}\end{align}

Actually, inequality (\ref{e107}) holds for any starting time, any
starting node $(i_1,i_2)$. Thus it holds that
$M_{({\mathcal{G}}(t),P_{\AS}(t))}=H_{({\mathcal{G}}(t)\times {\mathcal{G}}(t),Q(t),)}(\Theta) \leq 4Nt_1$.
This completes the proof.
\end{proof}

\subsection{Convergence analysis of $\AS$}

We are now in the position to characterize the convergence properties of $\AS$. The quantities $T_{\rm{con}}(x(0))$ and $T_{\Psi}$ for $\AS$ are defined in a similar way to those in Section~\ref{sec:AF}.

\begin{theorem} Let $x(0)\in {\mathbb{R}}^N$ and suppose $x(0)\notin
{\mathcal{W}}(x(0))$. Assume that ${\mathcal{G}}(t)$ be undirected and satisfies
Assumption~\ref{asm12}. Under $\AS$, almost any evolution $x(t)$ starting
from $x(0)$ reaches quantized average consensus. Furthermore,
${\mathbb{E}}[T_{\rm{con}}(x(0))]\leq\frac{1}{2}BJ(x(0))^2N^2(16N^7+1)$.\label{the6}
\end{theorem}

\begin{proof} Note that inequality \eqref{e106} also hold for $\AS$. Similar to Theorem~\ref{the2}, we have
${\mathbb{E}}[T_{\Psi}]=M_{({\mathcal{G}}(t),P_{\AS}(t))}$. Then,
the following estimate on ${\mathbb{E}}[T_{\rm{con}}(x(0))]$ holds:
\begin{align}{\mathbb{E}}[T_{\rm{con}}(x(0))]\leq \frac{N
J(x(0))^2}{8}M_{({\mathcal{G}}(t),P_{\AS}(t))}. \label{e192}\end{align}

Substituting the upper bound on $M_{({\mathcal{G}}(t),P_{\AS}(t))}$
in Proposition~\ref{pro2} into (\ref{e192}) and using $\log(\sqrt{2}N)\leq 2N$ gives the desired upper bound
on ${\mathbb{E}}[T_{\rm{con}}(x(0))]$ of $\AS$. The reminder of the proof on
the convergence to quantized average consensus is analogous to
Theorem~\ref{the4}, and thus omitted.
\end{proof}


\section{Discussion}

\subsection{Asynchronous distributed quantized averaging on random graphs}

Random graphs have been widely used to model real-world networks such as Internet,
transportation networks, communication networks, biological networks and social networks.
The Erd\H{o}s - R\'{e}nyi model ${\mathcal{G}}(N,p)$ is the most commonly studied
one, and constructed by randomly placing an edge between any two of
$N$ nodes with probability $p$.



For any given time, the probability that the (directed) edge $(i,j)$ is selected is\\
$p_0 := \frac{1}{N}\sum_{m=0}^{N-2}\frac{p}{m+1}C_{N-2}^m p^m (1-p)^{N-2-m}$, that is, node $i$ is active, the edge $(i,j)$ with other $m\in \{0,\cdots,N-2\}$ edges connecting node $i$ are placed, and the edge $(i,j)$ is selected by node $i$. To study the convergence properties of $\AF$ on ${\mathcal{G}}(N,p)$, it is equivalent to study $\AF$ on complete graphs with the transition matrix $P_{\AR}=(\hat{p}_{ij})\in {\mathbb{R}}^{N\times N}$ where $\hat{p}_{ij}=p_0$ and $\hat{p}_{ii}=1-(N-1)p_0$. The meeting time is denoted as $M_{({\mathcal{G}}(N,p),P_{\AR})}$. The probability that the two tokens meet for the first time at time $t$ is $2p_0$, that is, one of the tokens is chosen and simultaneously the edge between the two tokens is chosen. Hence, we have $M_{({\mathcal{G}}(N,p),P_{\AR})}=\sum_{\ell=1}^{=\infty}\ell2p_0(1-2p_0)^{\ell-1}=\frac{1}{2p_0}$.



Observe that $p_0=\frac{p}{N}\sum_{m=0}^{N-2}\frac{1}{m+1}C_{N-2}^m p^m
(1-p)^{N-2-m}\geq\frac{2p}{N(N-1)}\sum_{m=0}^{N-2}C_{N-2}^m p^m
(1-p)^{N-2-m}=\frac{2p}{N(N-1)}$. Like Theorem~\ref{the2}, we have
${\mathbb{E}}[T_{\rm{con}}(x(0))]\leq
\frac{NJ(x(0))^2}{8}{\mathbb{E}}[T_{\Psi}]=\frac{NJ(x(0))^2}{8}M_{({\mathcal{G}}(N,p),P_{\AR})}
=\frac{NJ(x(0))^2}{16p_0}\leq\frac{N^2(N-1)J(x(0))^2}{32p}.$

\subsection{Discussion on the bounds obtained}

Consider a fixed graph $L_N^m$ with $N$ vertices consists of a clique on $m$ vertices, including vertex $i$, and a path of length $N-m$ with one end connected to one vertex $k\neq i$ of the clique, and the other end of the path being $j$. It was shown in~\cite{Brightwell} that $H_{(L_N^{m_0}, P_{\SF})}$ is $O(N^3)$ where $m_0 = \lfloor\frac{2N+1}{3}\rfloor$. Let us consider the case that the algorithm $\AF$ is implemented on the graph $L_N^{m_0}$ and initial states $x_i(0) = 0$, $x_j(0) = 2$ and $x_k(0) = 1$ for all $k\neq i,j$. Observe that ${\mathbb{E}}[T_{\rm{con}}(x(0))] = M_{(L_N^{m_0}, P_{\AF})}$. From Proposition~\ref{pro1}, we have that ${\mathbb{E}}[T_{\rm{con}}(x(0))]$ is $O(N^4)$, that is one order less than the bound in Theorem~\ref{the2}.

Consider switching graphs ${\mathcal{G}}(t)$ where ${\mathcal{G}}(t)$ is the graph $L_N^{m_0}$ defined above when $t$ is a multiple of $B$; otherwise, all the vertices in ${\mathcal{G}}(t)$ are isolated. Random walks on ${\mathcal{G}}(t)$ can be viewed as time-scaled versions of those on $L_N^{m_0}$, that is, random walks on ${\mathcal{G}}(t)$ only make the movements when $t$ is a multiple of $B$. Let us consider the case that the algorithm $\AS$ is implemented on the graph $L_N^{m_0}$ and initial states $x_i(0) = 0$, $x_j(0) = 2$ and $x_k(0) = 1$ for all $k\neq i,j$. Following the same lines above, we have that the bound on ${\mathbb{E}}[T_{\rm{con}}(x(0))]$ is $O(BN^4)$ which is $N^4\log N$-order less than that in Theorem~\ref{the4}.

It can be directly computed that $H_{({\mathcal{G}}_{\rm com}, P_{\SF})}$ is $O(N^2)$ where ${\mathcal{G}}_{\rm com}$ is a complete graph with $N$ vertices. Following the same lines in Theorem~\ref{the2}, we have that ${\mathbb{E}}[T_{\rm{con}}(x(0))]$ is $O(N^3)$ when the algorithm $\AF$ is implemented on the graph ${\mathcal{G}}_{\rm com}$. It implies that the convergence of $\AF$ on ${\mathcal{G}}_{\rm com}$ is as fast as that on ${\mathcal{G}}(N,p)$ when $p$ is independent of $N$. This is consistent with the fact that the underlining graph of ${\mathcal{G}}(N,p)$ is ${\mathcal{G}}_{\rm com}$.

\bibliographystyle{plain}

\begin{thebibliography}{99}

\bibitem{Aysal} T.C. Aysal, M.J. Coates and M. Rabbat, {\em Distributed average
consensus using probabilistic quantization}, Proceedings of IEEE
Workshop on Statistical Signal Processing, Madison, USA, pp. 640 -
644, August 2007.
\bibitem{Bertsekas} D. P. Bertsekas and J. N. Tsitsiklis, {\em Parallel and distributed computation:
numerical methods}, Prentice Hall, 1989.
\bibitem{Blondel}V. D.
Blondel, J. M. Hendrickx, A. Olshevsky, and J. N. Tsitsiklis, {\em
Convergence in multiagent coordination, consensus, and flocking},
Proceedings of the Joint 44th IEEE Conference on Decision and
Control and European Control Conference (CDC-ECC'05), Seville,
Spain, pp. 2996 - 3000, December 2005.
\bibitem{Boyd} S. Boyd, A. Ghosh, B. Prabhakar,
and D. Shah, {\em Randomized gossip algorithms}, IEEE Transactions
on Information Theory, Special issue of IEEE Transactions on
Information Theory and IEEE ACM Transactions on Networking, vol. 52,
No. 6, pp. 2508 - 2530, June 2006.
\bibitem{Brightwell} G. Brightwell
and P. Winkler, {\em Maximum hitting time for random walks on
graphs}, Random Structures and Algorithms, vol. 1, No. 3, pp. 263 -
276, 1990.
\bibitem{Carli1} R. Carli, F. Bullo, and S. Zampieri, {\em Quantized average consensus via
dynamic coding/decoding schemes}, Procedings of IEEE Conference on Decision and Control, Cancun, Mexico, December 2008, To appear.
\bibitem{Carli2} R. Carli, F. Fagnani, P. Frasca, T. Taylor and S.
Zampieri, {\em Communication constraints in the state agreement
problem}, Automatica, to appear.
\bibitem{Coppersmith} D. Coppersmith, P. Tetali and P.
Winkler, {\em Collisions among random walks on a graph}, SIAM
Journal on Discrete Mathematics, vol. 6, No. 3, pp. 363 - 374, 1993.
\bibitem{Horn} R.A. Horn and C.R. Johnson, {\em Matrix analysis},
Cambridge, U.K.: Cambridge University Press, 1987.
\bibitem{Jadbabaie}A. Jadbabaie, J.
Lin, and A. S. Morse, {\em Coordination of groups of mobile
autonomous agents using nearest neighbor rules}, IEEE Transactions
on Automatic Control, vol. 48, No. 6, pp. 988 - 1001, June 2003.
\bibitem{Kashyap}A. Kashyap, T. Basar and R. Srikant, {\em Quantized consensus},
Automatica, vol. 43, pp. 1192 - 1203, July 2007.
\bibitem{Martinez} S. Mart\'{\i}nez, {\em Distributed
representation of spatial fields through adaptive interpolation
schemes}, Proceedings of the 2007 American Control Conference
(ACC'07), New York, USA, pp. 2750 - 2755, July 2007.
\bibitem{Nedic} A. Nedich, A. Olshevsky, A. Ozdaglar and J. N.
Tsitsiklis, {\em On distributed averaging algorithms and
quantization effects}, MIT LIDS Report 2778, November 2007,
https://netfiles.uiuc.edu/angelia/www/nedich.html.
\bibitem{Olfati1} R. Olfati-Saber, J. A. Fax, and R. M. Murray,
{\em Consensus and cooperation in networked multi-Agent systems},
Proceedings of the IEEE, vol. 95, No. 1, pp. 215-233, January 2007.
\bibitem{Olshevsky} A. Olshevsky and J.N. Tsitsiklis, {\em Convergence speed in
distributed consensus and averaging}, SIAM Journal on Control and Optimization, to appear.
\bibitem{Xiao} L. Xiao, S. Boyd and S. Lall, {\em A scheme for robust
distributed sensor fusion based on average consensus}, International
Conference on Information Processing in Sensor Networks, pp. 63 -
70, Los Angeles, 2005.
\bibitem{Zhu} M. Zhu and S. Mart\'{\i}nez, {\em On the convergence time of distributed
quantized averaging algorithms}, Proceedings of the 47th IEEE Conference on Decision and Control,
Canc\'{u}n, Mexico, to appear 2008
\end{thebibliography}


\end{document}